\documentclass{article}
\usepackage{amsfonts,amsmath,amssymb,amsbsy}
\usepackage{graphicx}
\usepackage{epsfig}
\usepackage{theorem}
\newtheorem{theorem}{Theorem}
\newtheorem{lemma}{Lemma}
\newtheorem{cor}{Corollary}
\newtheorem{defin}{Definition}

\begin{document}

\begin{titlepage}


\title{\textbf{Fractal analysis of hyperbolic and nonhyperbolic fixed points and singularities of dynamical systems in $\mathbb{R}^{n}$}}

\author{Lana Horvat Dmitrovi\' c\\
University of Zagreb Faculty of Electrical Engineering and Computing,\\ 
Unska 3, 10000 Zagreb, Croatia,\\
lana.horvat@fer.hr}

\maketitle

\begin{abstract}
The main purpose of this article is to study box dimension of orbits near hyperbolic and nonhyperbolic fixed points of discrete dynamical systems in $\mathbb{R}^{n}$. We generalize the known results for one-dimensional systems, that is, the orbits near the hyperbolic fixed point in one-dimensional discrete dynamical system has the box dimension equal to zero and the orbits near nonhyperbolic fixed point has positive box dimension. In the process of studying box dimensions, we use the stable, unstable and center manifolds and appropriate system restrictions. The main result is that the box dimension of orbit equals zero near stable and unstable hyperbolic fixed points and on the stable and unstable manifolds. The main results in the nonhyperbolic case is that box dimension is determined by the box dimension on the center manifold. We also introduce the projective box dimension and use it as a sufficient condition for nonhyperbolicity. At the end, all the results for discrete systems can be applied to continuous systems by using the unit time map so we apply it to the hyperbolic and nonhyperbolic singularities of continuous dynamical systems in $\mathbb{R}^{n}$.
\end{abstract}

\end{titlepage} 

\def\b{\beta}
\def\g{\gamma}
\def\d{\delta}
\def\l{\lambda}
\def\o{\omega}
\def\ty{\infty}
\def\e{\varepsilon}
\def\f{\varphi}
\def\:{{\penalty10000\hbox{\kern1mm\rm:\kern1mm}\penalty10000}}
\def\st{\subset}
\def\stq{\subseteq}
\def\q{\quad}
\def\M{{\cal M}}
\def\cal{\mathcal}
\def\eR{\mathbb{R}}
\def\eN{\mathbb{N}}
\def\Ze{\mathbb{Z}}
\def\Qu{\mathbb{Q}}
\def\Ce{\mathbb{C}}
\def\ov#1{\overline{#1}}
\def\D{\Delta}
\def\O{\Omega}

\def\bg{\begin}
\def\eq{equation}
\def\bgeq{\bg{\eq}}
\def\endeq{\end{\eq}}
\def\bgeqnn{\bg{eqnarray*}}
\def\endeqnn{\end{eqnarray*}}
\def\bgeqn{\bg{eqnarray}}
\def\endeqn{\end{eqnarray}}

\newcount\remarkbroj \remarkbroj=0
\def\remark{\advance\remarkbroj by1 \smallskip{Remark\ \the\remarkbroj.}\enspace\ignorespaces\,}

\pagestyle{myheadings}{\markright{Characterisation of hyperbolic/nonhyperbolic fixed point}}

\section{Introduction}

The main motivation for this study is that in the article \cite{neveda} was shown that the box dimension of the orbits of one-dimensional discrete system near the hyperbolic fixed points is zero, while in the article \cite{laho} the result that the box dimension near the nonhyperbolic fixed point is stricly positive is presented. Therefore, we would like to get the analogous result for the hyperbolic and nonhyperbolic fixed points of the discrete systems in $\mathbb{R}^{n}$. Of course, this is a first step in doing the fractal analysis of bifurcations in $\mathbb{R}^{n}$. Article \cite{laho1} started with studying the bifurcations of discrete dynamical systems in $\mathbb{R}^{2}$. It was shown that in the case of one multiplier on the unit circle, the box dimension of the orbits on the center manifold is connected to the bifurcation. That is, saddle-node and transcritical bifurcations has box dimension $\frac{1}{2}$ and pitchfork and period-doubling bifurcations are connected to the box dimension $\frac{2}{3}$. While in the case of Neimark-Sacker, the box dimension in the rational case is $\frac{2}{3}$ and in the irrational case is $\frac{4}{3}$. 

The aim of this paper  is the generalisation of results about box dimension near the hyperbolic and non-hyperbolic fixed points in higher dimensions. In the process of determining box dimension, the restriction of systems to stable, unstable and center manifolds, Lipschitz property of box dimension and the notion of characteristic box dimension are used. This specific change in box dimension of an orbit at the moment of bifurcation has already been explored for some bifurcations in one and two dimensions. This paper further explores this connection of box dimension as fractal property to some bifurcations in higher dimensions, such as bifurcations with onedimensional center manifold. Moreover, we need to mention the article which introduced the fractal analysis of dynamical systems and bifurcations, that is, $\cite{zuzu}$. In that article the result of box dimension for Hopf bifurcation is proved and the change of box dimension at the bifurcation point has been noticed for the first time.

Now we recall the notions of box dimension and Minkowski content. For further details see e.g. \cite{fa}, \cite{zu2}.
For bounded set $A\st\eR^N$ {\it Minkowski neighbourhood} of radius $\e$ around $A$ is a $\e$-neighbourhood of $A$, that is
$A_\e=\{y\in\eR^N\:d(y,A)<\e\}$. For $s\ge0$,\textit{the lower and upper $s$-dimensional Minkowski contents of $A$} are defined by
$$\M_*^s(A):=\liminf_{\e\to0}\frac{|A_\e|}{\e^{N-s}},\,\,\,\,\M^{*s}(A):=\limsup_{\e\to0}\frac{|A_\e|}{\e^{N-s}}.$$
Then \textit{the lower and upper box dimension} are defined by 
$$\underline\dim_BA=\inf\{s>0:\M_*^s(A)=0\}, \,\,\,\,\ov\dim_BA=\inf\{s>0:\M^{*s}(A)=0\}.$$
If $\underline\dim_BA=\ov\dim_BA$ we denote it by $\dim_BA$.
If there exists $d\ge0$ such that\ $0<\M_*^d(A)\le \M^{*d}(A)<\ty,$ then we say that set $A$ is \textit{Minkowski nondegenerate}. Clearly, then $d=\dim_B A$. If $|A_\e|\simeq \e^{s}$ for $\e$ small, then $A$ is Minkowski nondegenerate set and 
$\dim_B A=N-s$. 

We say that $\mathbf{f}:A\rightarrow B$, where $A, B$ are open subsets of $\mathbb{R}^{n}$, is a \textit{bilipschitz map} if there exist positive constants $C_1$ and $C_2$ such that
$$C_1\left\|\mathbf{x}-\mathbf{y}\right\|\leq\left\|\mathbf{f}(\mathbf{x})-\mathbf{f}(\mathbf{y})\right\|\leq C_2 \left\|\mathbf{x}-\mathbf{y}\right\|$$
for every $\mathbf{x}, \mathbf{y}\in A$. If $\mathbf{f}$ is a bilipschitz mapping, than $\dim_{B}A=\dim_{B}\mathbf{f}(A).$



The remainder of this paper is organized as follows. In the second section we present the main result about the box dimension of the orbit of discrete dynamical system around the hyperbolic fixed point. In Section 3 we study the box dimension near the nonhyperbolic fixed point. In Section 4 we apply the result to continuous dynamical systems with hyperbolic and nonhyperbolic singularities via unit-time map. 


\section{Box dimension of hyperbolic fixed point in $\mathbb{R}^{n}$}

We consider the box dimension of an orbit in the neighbourhood of the hyperbolic fixed point.
We will see that the box dimension is connected to the hyperbolicity and nonhyperbolicity of a fixed point. 
It was shown in the article \cite{neveda} that the box dimension of hyperbolic fixed point of one-dimensional discrete dynamical system is zero. We will now consider $n$-dimensional systems.

We consider the discrete dynamical system
\begin{eqnarray} \label{s2}
\mathbf{x}_{n+1}=\mathbf{F}(\mathbf{x}_{n})
\end{eqnarray}
with $\mathbf{F}:\mathbb{R}^{n}\rightarrow \mathbb{R}^{n}$ and $\mathbf{x}_{n}\in \mathbb{R}^{n}$.

Since we look at the local behaviour in the neighborhood of a fixed point, we presume that $\mathbf{F}$ is smooth enough, that is, of class $C^{r}$ on that neighbourdood. The fixed point of a map $\mathbf{F}$ is $\mathbf{x}_0\in \mathbb{R}^{n}$ such that $\mathbf{F}(\mathbf{x}_0)=\mathbf{x_0}$. If  $A=D\mathbf{F}(\mathbf{x_0})$ is a Jacobi matrix in the point $\mathbf{x_0}$, then its eigenvalues $\lambda_{i}$, $i=1,\ldots,n$ are called the multipliers of a fixed point. Fixed point $\mathbf{x_0}$ is a hyperbolic fixed point if $\left|\lambda_i\right|\neq 1$ for $i=1,\ldots,n$. 
We also know that the hyperbolic fixed point is stable if all the multiplicators are inside the unit circle.
We will recall the definition of exponential stability (see \cite{io}).
\begin{defin}
Fixed point $\mathbf{x_0}=0$ of a map $\mathbf{F}$ is \textbf{exponentially stable} if there exists the neighbourhood $V$, constants $\gamma>0$ i $K\in(0,1)$ such that $\mathbf{x_0}=0$ is Lyapunov stable and for every $\mathbf{x}\in V$ is $\left\|\mathbf{F}^{n}(\mathbf{0})\right\|\leq \gamma K^{n}$, for all $n$.
\end{defin}

\begin{theorem} (\cite{io})\\
Let $\mathbf{F}:\mathbb{R}^{n}\rightarrow \mathbb{R}^{n}$ is of class $C^1$ and $\mathbf{x_0}=0$ is a fixed point of $\mathbf{F}$, and let $A=D\mathbf{F}(\mathbf{0})$ be a Jacobi matrix of a map $\mathbf{F}$ in $\mathbf{x_0}$. If the spectrum of $A$ lies inside the open unit circle, then the fixed point $\mathbf{x_0}$ is exponentially stable.
\end{theorem}

Now we can prove the theorem about the box dimension of every orbit in the neighbourhood of the stable hyperbolic fixed point in $\mathbb{R}^{n}$.

\begin{theorem} \textbf{Stable hyperbolic fixed point in $\mathbb{R}^{n}$}\\
Let $\mathbf{F}:I\rightarrow \mathbb{R}^{n}$, $I=B_{r}(\mathbf{x_0})\subset\mathbb{R}^{n}$ be a map of class $C^1$ such that $\mathbf{F}(\mathbf{x_0})=\mathbf{x_0}$ and all the eigenvalues of $D\mathbf{F}(\mathbf{x_0})$ lies inside the unit circle, that is, $\left|\lambda_{i}\right|<1$, for $i=1,\ldots,n$. Let $S(\mathbf{x_1})=(\mathbf{x}_{n})_{n\in\mathbb{N}}$ be an orbit of the discrete dynamical system generated by $\mathbf{x}_{n+1}=\mathbf{F}(\mathbf{x}_{n})$ with $\mathbf{x_1}\in B_{r}(\mathbf{x_0})$.  Then there exists $r_1<r$ such that for every $\mathbf{x_1}\in B_{r_1}(\mathbf{x_0})$ we have $\dim_{B} S(\mathbf{x_1})=0$.
\end{theorem}

\textbf{Proof.}\\
Without loss of generality we may assume that the fixed point is in the origin  $\mathbf{x_0}=0$. From the assumptions and Theorem 1 it follows that the origin is exponentially stable, that is, there exists the neighbourhood $V$, constants $\gamma>0$ and $k\in(0,1)$ such that for every $n$ it holds
$$\left\|\mathbf{F}^{n}(\mathbf{x})\right\|\leq \gamma k^{n}.$$
Also we get $\left\|\mathbf{F}^{n}(\mathbf{x})\right\|\leq\left\|\mathbf{x}\right\|$ for every $\mathbf{x}\in V$, what means that the sequence $S(\mathbf{x_1})$ is decreasing by the norm on $V$.\\
Now we consider the Minkowski neighbourhood of a set $S=S(\mathbf{x_1})$, where $\mathbf{x_1}\in V$. For every fixed $\varepsilon>0$ the inequality
$$\left\|\mathbf{x}_{n}\right\|=r_{n}\leq \gamma k^{n}$$ is satisfied for  $n\geq n_0(\varepsilon):=\left\lceil \frac{\log \frac{2\varepsilon}{\gamma}}{\log k}\right\rceil$. Now for $\varepsilon$-Minkowski neighbourhood of the set $S=S(\mathbf{x_1})$ the next inequality holds
\begin{eqnarray}
\left|S_{\varepsilon}\right|&\leq& c_{n}\varepsilon^{n} n_0(\varepsilon) + c_{n}\varepsilon^{n} \leq c_{n}\varepsilon^{n}(\frac{\log \frac{2\varepsilon}{\gamma}}{\log k}+1)+c_{n}\varepsilon^{n}\\
&=&(2-\frac{\log \gamma}{\log k}+\frac{\log 2}{\log k})c_{n}\varepsilon^{n} + \frac{\log \varepsilon}{\log k}c_{n}\varepsilon^{n} \nonumber
 \end{eqnarray}
So we have
\begin{eqnarray}
\frac{\left|S_{\varepsilon}\right|}{\varepsilon^{n-s}}\leq c_{n} (2-\frac{\log \gamma}{\log k}+\frac{\log 2}{\log k})\varepsilon^{s}+\frac{c_{n}}{\log k}\varepsilon^{s}\log \varepsilon.
\end{eqnarray}
It follows that $\mathcal{M}^{*s}(S)=0$ for every $s\in(0,1]$ and the claim is proven, that is, $\dim_{B}S=0$.$\blacksquare$\\

\begin{cor} \textbf{Unstable hyperbolic fixed point in $\mathbb{R}^{n}$}\\
Let $\mathbf{F}:I\rightarrow \mathbb{R}^{n}$, $I=B_{r}(\mathbf{x_0})\subset\mathbb{R}^{n}$ be a map of class $C^1$ such that $\mathbf{F}(\mathbf{x_0})=\mathbf{x_0}$ and all the eigenvalues of $D\mathbf{F}(\mathbf{x_0})$ lies outside the unit circle, that is, $\left|\lambda_{i}\right|>1$, for $i=1,\ldots,n$. Let $S(\mathbf{x_1})=(\mathbf{x}_{n})_{n\in\mathbb{N}}$ be an orbit of the discrete dynamical system generated by $\mathbf{x}_{n+1}=\mathbf{F}^{-1}(\mathbf{x}_{n})$ with $\mathbf{x_1}\in B_{r}(\mathbf{x_0})$.  Then there exists $r_1<r$ such that for every $\mathbf{x_1}\in B_{r_1}(\mathbf{x_0})$ we have $\dim_{B} S(\mathbf{x_1})=0$.
\end{cor}

\textbf{Proof.}\\
We know that if the operator $D\mathbf{F}(\mathbf{x_0})$ doesn't have zero eigenvalue, then it is a regular operator at $\mathbf{x_0}$.
Then by the Inverse Function Theorem, we have that the map $F$ of class $C^1$ on $B_{r}(\mathbf{x_0})$ has a inverse map $F^{-1}$ 
on some neighbourhood of $\mathbf{x_0}=\mathbf{F}(\mathbf{x_0})$ of the same class $C^1$, and it holds that $(DF^{-1})(F(\mathbf{x_0}))=DF(\mathbf{x_0})^{-1}$. Since the eigenvalues of the inverse operator are $\mu_{i}=\frac{1}{\lambda_{i}}$, then the eigenvalue of $DF^{-1}(\mathbf{x_0})$ lies inside the unit circle, and we can apply Theorem 2. $\blacksquare$\\

\textbf{Example 1.} We present some cases of hyperbolic fixed points at the origin in $\mathbb{R}^2$ such as stable node and stable focus. At Figure 1a is the phase portrait of the system 
\begin{eqnarray} 
x_{n+1}&=&0.8\, x_{n} \nonumber\\
y_{n+1}&=&0.7\,y_{n} \nonumber
\end{eqnarray}
which has a stable node at the origin with two real positive eigenvalues inside the unit circle. At Figure 1b is the phase portrait of the system
\begin{eqnarray} 
x_{n+1}=&0.8\, x_{n} &+ 0.2\,y_{n}\nonumber\\
y_{n+1}=&-0.2\,x_{n} &+ 0.8\,y_{n}\nonumber
\end{eqnarray}
which has a stable focus with a pair of complex conjugate multipliers inside the unit circle. Therefore, by Theorem 1 in both cases we have $\dim_{B}S(x_1,y_1)=0$ for every orbit near the origin.\\

\begin{center}
\includegraphics[width=5cm]{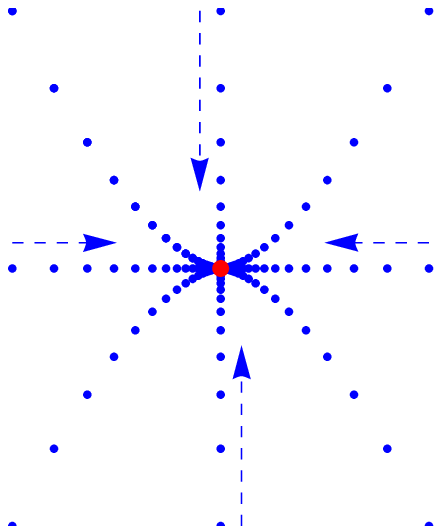} \hskip 1cm\includegraphics[width=5cm]{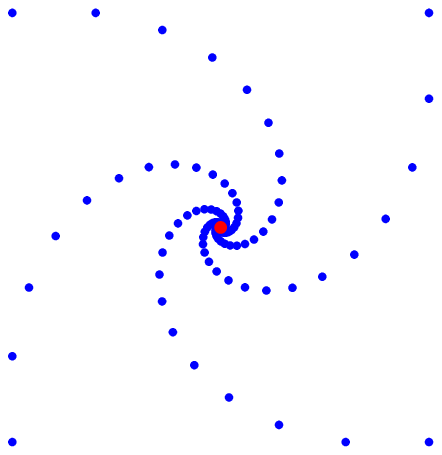}\\ 
 \textbf{Figure 1a} - hyperbolic node \hskip 2cm  \textbf{Figure 1b} - hyperbolic focus
\end{center}

\subsection{Stable and unstable manifolds in $\mathbb{R}^2$}

In this section we will consider the stable and unstable manifolds in $\mathbb{R}^{2}$. We need the following lemma from article \cite{zuzu3}.

\begin{lemma} \cite{zuzu3}\\
Let $g:\mathbb{R}^{N-1}\rightarrow \mathbb{R}$ be a Lipchitz map, $N\geq 2$, and we define $F:\mathbb{R}^{N}\rightarrow\mathbb{R}^{N}$ with $F(X,Z):=(X,Z+g(X))$, where $X\in\mathbb{R}^{N-1}$ and $Z\in\mathbb{R}$. Then $F$ is a bilipschitz map and measure preserving, that is, for every measurable set $E\subset\mathbb{R}^{N}$ of limited measure it holds $\left|F(E)\right|=\left|E\right|$. Furthermore, for every limited set $A\subset\mathbb{R}^{N}$ we have
$$\overline{\dim}_{B}F(A)=\overline{\dim}_{B}A\,\,\, \rm{and}\,\,\, \underline{\dim}_{B}F(A)=\underline{\dim}_{B}A.$$
Set $A$ is a nondegenerate if and only if $F(A)$ is a nondegenerate.
\end{lemma}

If we apply this lemma to the plane ($N=2$), then for the Lipschitz map $g:\mathbb{R}\rightarrow\mathbb{R}$ the map
$F:\mathbb{R}^2\rightarrow \mathbb{R}^2$ defined by
\begin{equation} \label{s3}
F(x,z)=(x,z+g(x))
\end{equation}
 is a bilipschitz map.\\

We consider two-dimensional system 
\begin{equation} \label{jed1}
\left(x,y\right)\longmapsto F(x,y),\,\,F:\mathbb{R}^2\rightarrow\mathbb{R}^2
\end{equation} 
with hyperbolic fixed point at $(x,y)=(0,0)$. Now by using Lemma 1 we can prove the theorem with box dimension on the stable and unstable manifolds.

\begin{theorem} \textbf{Box dimension on the stable manifold}\\
Let the restriction of the system $(\ref{jed1})$ with $\left|\lambda_1\right|<1$ on the stable manifold $y=h(x)$ is given by
$$x\mapsto G(x)=\lambda_1 x + f(x,h(x))$$
where function $f$ has Taylor expansion starting with at least quadratic terms.
Let $S(x_1,y_1)=(x_{n},y_{n})$ be an orbit of the system ($\ref{jed1}$) on the stable manifold in the form 
\begin{eqnarray} \label{jed2}
x_{n+1}&=&\lambda_1 x_{n} + f(x_{n},h(x_{n}))\\
y_{n+1}&=&h(x_{n+1}) \nonumber
\end{eqnarray}
with initial point $(x_1,y_1)$. 
Then there exists $r>0$ such that for $\left\|(x_1,y_1)\right\|<r$ we have $$\dim_{B}S(x_1,y_1)=0.$$
\end{theorem}

\textbf{Proof.}\\
First we define the set $A=\{(x,y):x\in A(x_1),y=0\}$, where $A(x_1)=(x_{n})_{n\in\mathbb{N}}$ is one-dimensional discrete dynamical system generated by $x_{n+1}=G(x_{n})=\lambda_1 x_{n} + f(x_{n},h(x_{n}))$ and $x_1\in(0,r)$. Then we act with the map $F(x,z)=(x,z+h(x))$ on the set $A$ and get $$F(A)=F(x,0)=(x,0+h(x))=(x,h(x)).$$ 
So the image of the set $A$ under map $F$ is associated to the system $(\ref{jed2})$ on the center manifold $y=h(x)$. In other words, the map $F$ 
associates the projection of an orbit on $x$-axis with appropriate orbit on the center manifold. 
Since the map $h$ is of class $C^{r}$ on some neighbourhood small enough $\left|x\right|<\delta$ and $h'(0)=0$, then $h'$ is limited on that neighbourhood, and by the Mean Value Theorem we have
$$\left|h(x_2)-h(x_1)\right|=\left|h'(x^{*})\right|\left|x_2-x_1\right|\leq M\left|x_2-x_1\right|$$
for some $x^{*}\in(x_1,x_2)$. Therefore $h$ is a Lipschitz map for $\left|x\right|<\delta$.
Now it follows from Lemma that $F$ is a bilipschitz map, and 
$$\dim_{B} F(A)=\dim_{B}A.$$
Since $h$ has a stable fixed point at $x_0=0$, from \cite{neveda}, Theorem 3.1 follows that there exists $r>0$ such that for the sequence $A(x_1)=(x_{n})_{n\in\mathbb{N}}$ defined by $x_{n+1}=h(x_{n})$, $\left|x_1\right|<r$ we have
$$\dim_{B} A(x_1)=0.$$
We see now that $\dim_{B}A=\dim_{B}A(x_1)=0$.
Notice that $S(x_1,y_1)=F(A)$, so it follows that for $|x_1|<r$ is
$$\dim_{B}S(x_1,y_1)=\dim_{B} F(A)=\dim_{B}A=0.\,\,\,\,\blacksquare$$

\textbf{Remark 1.} Analogous result holds for the inverse orbit on one-dimensional unstable manifold.\\

\textbf{Example 2.} \textbf{Hyperbolic saddle in $\mathbb{R}^{2}$}\\
We will show how the above procedure works in the case $0<\lambda_1<1<\lambda_2$. The other cases with negative multipliers will be analogous 
but with alternating orbits.

Well, we consider the map $\mathbf{F}:\mathbb{R}^2\rightarrow \mathbb{R}^2$ defined by
\begin{eqnarray} \label{sed1}
x&\mapsto& \lambda_1 x + a_1x^2 + a_2 xy + b_2 y^2 + \mathcal{O}(\left|x\right|+\left|y\right|)^{3}\nonumber \\
y&\mapsto& \lambda_2 y + b_1 x^2 + b_2 xy + b_3 y^2 + \mathcal{O}(\left|x\right|+\left|y\right|)^{3},
\end{eqnarray}
where $\lambda_1<1<\lambda_2$. Stable manifold is given by $W^{s}=\{(x,y): y=\phi_1(x)\}$, and unstable by $W^{u}=\{(x,y):x=\phi_2(y)\}$. We write
\begin{eqnarray}
\phi_1(x)&=&\alpha_2 x^2 + \alpha_3 x^3 + \mathcal{O}(\left|x\right|^4)\\
\phi_2(y)&=&\beta_2 y^2 + \beta_3 y^3 + \mathcal{O}(\left|y\right|^4),
\end{eqnarray}
and get
$$\alpha_2=\frac{b_1}{\lambda_1^2-\lambda_2},\,\,\,\alpha_3=\frac{b_1(b_2-2a_1\lambda_1)}{(\lambda_1^2-\lambda_2)(\lambda_1^3-\lambda_2)},\,\,\,\beta_2=\frac{a_3}{\lambda_2^2-\lambda_1},\,\,\,\beta_3=\frac{a_3(a_2-2b_3\lambda_2)}{(\lambda_2^2-\lambda_1)(\lambda_2^3-\lambda_1)}.$$
Now we see that the restriction of the system ($\ref{sed1}$) on the stable manifold $W^{s}$ is
\begin{eqnarray} \label{res1}
x_{n+1}&=& \lambda_1 x_{n} + a_1x_{n}^2 + a_2 \alpha_2 x_{n}^3 + \mathcal{O}(\left|x_{n}\right|^{4})\nonumber \\
y_{n+1}&=& \phi_1(x_{n+1}),
\end{eqnarray}
while the restriction of ($\ref{sed1}$) on the unstable manifold $W^{u}$ is given by
\begin{eqnarray} \label{res2}
x_{n+1}&=&\phi_2(y_{n+1})\nonumber \\
y_{n+1}&=& \lambda_2 y_{n} + b_3 y_{n}^2 + b_2 \beta_2 y_{n}^3 +  \mathcal{O}(\left|y_{n}\right|^{4}).
\end{eqnarray}

Since the maps $\phi_1$ and $\phi_2$ are of class $C^{r}$, then they will be Lipschitz maps on some neighbourhood of a fixed point $(0,0)$. By using Lemma 1 we get that the orthogonal projection of the restriction ($\ref{res1}$) on the $x$-axis is a bilipschitz map. Therefore the box dimension of the orbit on $W^{s}$ is equal to the box dimension of the orbit of one-dimensional discrete dynamical system
$$x_{n+1}=\lambda_1 x_{n} + a_1x_{n}^2 + a_2 \alpha_2 x_{n}^3 + \mathcal{O}(\left|x_{n}\right|^{4}).$$
Analogously for the orbit on $W^{u}$. At this moment, we will use the fact that in the neighbourhood of hyperbolic fixed point of one-dimensional systems the box dimension of all orbits equals zero. 

So, for the twodimensional sequence $S(x_1,y_1)=(x_{n},y_{n})$ defined by ($\ref{res1}$) there exists $r>0$ such that for $\left|x_1\right|<r$ we have
$\dim_{B} S(x_1,y_1)=0$. Analogously we get the result for ($\ref{res2}$). We established that the box dimension of the orbit on the stable and unstable manifold in the neighbourhood of hyperbolic saddle $(0,0)$ are zero.

At Figure 2 we see the phase portrait of a system $x_{n+1}=1.2\,x_{n},\,\,y_{n+1}=0.7\,y_{n}$ with stable manifold $x=0$ and unstable $y=0$. Box dimensions of orbits on the manifolds are zero.

\begin{center}
\includegraphics[width=5cm]{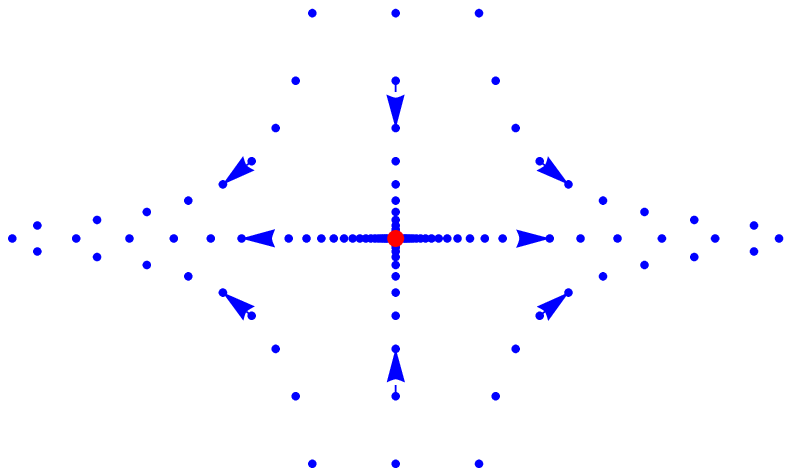} \\ 
\textbf{Figure 2} - hyperbolic saddle
\end{center}

\subsection{Stable manifolds in $\mathbb{R}^{n}$}

Now we consider the hyperbolic fixed point with stable and unstable manifolds in $\mathbb{R}^{n}$. 

\begin{theorem} \textbf{Local stable and unstable manifold}, \cite{guho}\\
Let $G:\mathbb{R}^{n}\rightarrow \mathbb{R}^{n}$ be a $C^{r}$-diffeomorphism with a hyperbolic fixed point $x_0$, with $n_0=0$, $n_{+}+n_{-}=n$. Then the intersections of $W^{s}(\mathbf{x_0})$ and $W^{u}(\mathbf{x_0})$ with a sufficiently small neighbourhood of $\mathbf{x_0}$ contain $C^{r}$ smooth manifolds $W^{s}_{loc}(\mathbf{x_0})$ and $W^{u}_{loc}(\mathbf{x_0})$ of dimension $n_{-}$ and $n_{+}$, respectively. Moreover, $W^{s}_{loc}(\mathbf{x_0})$ ($W^{u}_{loc}(\mathbf{x_0})$) is tangent at $\mathbf{x_0}$ to $T^{s}$($T^{u}$), where $T^{s}$($T^{u}$) is  the generalized eigenspace corresponding to the union of all eigenvalues of $A$ with $\left|\lambda\right|<1$ ($\left|\lambda\right|>1$).
\end{theorem}
 
\textit{Stable} (\textit{unstable}) \textit{invariant set} for a fixed point $\mathbf{x_0}$ is defined by
$$W^{s}(\mathbf{x_0})=\{\mathbf{x}:\mathbf{f}^{k}(\mathbf{x})\rightarrow \mathbf{x_0},\,k\rightarrow \infty\},$$
$$W^{u}(\mathbf{x_0})=\{\mathbf{x}:\mathbf{f}^{k}(\mathbf{x})\rightarrow \mathbf{x_0},\,k\rightarrow -\infty\}.$$

\textbf{Remark 2.} Notice that if $n_{+}=n$ then the point is unstable, and if $n_{-}=n$ is stable, so we can apply Theorem 1 and Corollary 1.

In order to study the stable manifolds in $\mathbb{R}^{n}$, we need to generalise the result from Lemma 1.

\begin{lemma}
Let $G:\mathbb{R}^{m}\rightarrow \mathbb{R}^{N-m}$ be a map of class $C^1$, $N\geq 2$, $m\leq N$ and we define $F:\mathbb{R}^{N}\rightarrow\mathbb{R}^{N}$ with $F(X,Z):=(X,Z+G(X))$, where $X\in\mathbb{R}^{m}$ i $Z\in\mathbb{R}^{N-m}$, $G(X)=(g_1(X),\ldots,g_{N-m}(X))$. Then $F$ is a bilipschitz map and  for every limited set $A\subset\mathbb{R}^{N}$ we have
$$\overline{\dim}_{B}F(A) = \overline{\dim}_{B}A\,\,\, \rm{and}\,\,\, \underline{\dim}_{B}F(A) = \underline{\dim}_{B}A.$$
\end{lemma}

\textbf{Proof.} Since $G$ is $C^1$-map, it follows that $F$ is also a $C^1$-map and we know that on the limited set $F$ is a Lipschitz map. It is easy to see that $DF(0)$ is regular operator, so there exists $F^{-1}$ which is also of class $C^1$ and Lipschitz map. $\square$

\begin{lemma} 
Let $P:\mathbb{R}^{N}\rightarrow\mathbb{R}^{N}$ be a m-dimensional orthogonal projection defined with $P(X,Z):=(X,0)$, where $X\in\mathbb{R}^{m}$, and $N> m\geq 1$. Then $P$ is a Lipschitz map and  for every limited set $A\subset\mathbb{R}^{N}$ we have
$$\overline{\dim}_{B}P(A)\leq \min\left\{m,\overline{\dim}_{B}A\right\}\,\,\, \rm{and}\,\,\, \underline{\dim}_{B}P(A)\leq \min\left\{m,\underline{\dim}_{B}A\right\}.$$
\end{lemma}

\textbf{Remark 3.} We will also use this lemma in the case of nonhyperbolic points. That is, if the projection of orbit on any coordinate line has stricly positive box dimension, that the box dimension of orbit is also stricly positive.

Now by using the above lemma we can prove the theorem about box dimension on the stable manifold.

\begin{theorem} \textbf{Box dimension on the stable manifold}\\
Let we have a system $\mathbf{x}\mapsto \mathbf{F}(\mathbf{x})$, $\mathbf{F}:\mathbb{R}^{N}\rightarrow \mathbb{R}^{N}$. Let $\mathbf{x_0}=0$ be a hyperbolic fixed point with $m$ multipliers 
inside the unit circle, and $N-m$ multipliers outside the unit circle. Let the stable manifold be given by $Y=G(X)$, where  $G:\mathbb{R}^{m}\rightarrow \mathbb{R}^{N-m}$ is of class $C^{r}$.
Let $S(X_1,Y_1)=\{(X_{n},Y_{n})\}_{n\in\mathbb{N}}$ be an orbit of the systems on the stable manifold in the form 
\begin{eqnarray}
X_{n+1}&=&\Lambda X_{n} + F_2(X_{n},G(X_{n}))\\
Y_{n+1}&=&G(X_{n+1}) \nonumber
\end{eqnarray}
where initial point $(X_1,Y_1)$. 
Then there exists $r>0$ such that for $\left\|(X_1,Y_1)\right\|<r$ we have $$\dim_{B}S (X_1,Y_1)=0.$$
\end{theorem}

\textbf{Proof.}\\
For $X\in\mathbb{R}^{m}$ and $Y\in\mathbb{R}^{N-m}$ we define a set $A=\{(X,Y)\in\mathbb{R}^{m}\times\mathbb{R}^{N-m} : X\in A(X_1), Y=0\}$, where $A(x_1)=(X_{n})_{n\in\mathbb{N}}$ is $m$-dimensional discrete dynamical system generated by $X_{n+1}=h(X_{n})=\Lambda X_{n} + F_2(X_{n},G(X_{n}))$ and $\left\|(X_1,0)\right\|<r$. Then we act with the map $Q$ on the set $A$ that is with  $Q(X,Z)=(X,G(X))$ and get $$Q(A)=Q(X,0)=(X,G(X)).$$ 
So the image of the set $A$ under map $Q$ is associated to the system $(12)$ on the center manifold $Y=G(X)$. In other words, the map $Q$ 
associate the projection of an orbit on $m$-dimensional coordinate plane with appropriate orbit on the center manifold. 
Since the map $G$ is of class $C^{r}$, $r\geq 1$ on some neighbourhood small enough $\left\|X\right\|<\delta$, then $\left\|D G(0)\right\|\leq M$, $M>0$ ($DG$ is limited operator on that neighbourhood), and we have
$$\left\|G(X_2)-G(X_1)\right\|=\left\|D G(X^{*})\right\|\left\|X_2-X_1\right\|\leq M\left\|X_2-X_1\right\|$$
for some $X^{*}$ on the line between $X_1$ and $X_2$. Therefore $G$ is Lipschitz map for $\left\|X\right\|<\delta$.
Now it follows from Lemma 2 that $F$ is a Lipschitz map, and 
$$\dim_{B} F(A)\leq \dim_{B}A.$$
Since $G$ is a $C^1$ map with all the eigenvalues of $DG(0)$ inside the unit circle, so from Theorem 1 follows that there exists $r>0$ such that for the sequence $A(X_1)=(X_{n})_{n\in\mathbb{N}}$ defined by $X_{n+1}=G(X_{n})$, $\left\|X_1\right\|<r$ we have
$$\dim_{B} A(X_1)=0.$$
We see now that $\dim_{B}A=\dim_{B}A(x_1)=0$.
Notice that $S(X_1,Y_1)=Q(A)$, so it follows that for $\left\|(X_1,Y_1)\right\|<r$ is
$$\dim_{B}S(X_1,Y_1)=\dim_{B} Q(A)\leq \dim_{B}A=0.\,\,\,\,\blacksquare$$


\begin{defin}(\textbf{Projective box dimension})
Let we have a system $X\mapsto F(X)$, $F:\mathbb{R}^{N}\rightarrow \mathbb{R}^{N}$, and let $S(X_1)=(X_{n})_{n\in\mathbb{N}}$ be an orbit of the system  with initial point $X_1$. Let $P_{k}:\mathbb{R}^{N}\rightarrow\mathbb{R}^{N}$, $k=1,\ldots,N$ be onedimensional projections defined with $P(x_1,\ldots,x_{N}):=(0,\ldots,0,x_{k},0\ldots,0)$. Then we define projective box dimensions of the orbit $S(X_1)$ with $$\dim_{P_{k}}S(X_1)=\dim_{B}P_{k}(S(X_1)),\,\,k=1,\ldots,N.$$
\end{defin}

\begin{cor}(\textbf{Projective box dimensions of stable hyperbolic fixed point})\\
Let we have a system $\mathbf{x}\mapsto \mathbf{F}(\mathbf{x})$, $\mathbf{F}:\mathbb{R}^{N}\rightarrow \mathbb{R}^{N}$. Let $\mathbf{x_0}=0$ be a hyperbolic fixed point with $N$ multipliers inside the unit circle.
Let $S(X_1)=(X_{n})_{n\in\mathbb{N}}$ be an orbit of the system $X\mapsto F(X)$ with initial point $X_1$. 
Let $P_{k}:\mathbb{R}^{N}\rightarrow\mathbb{R}^{N}$, $k=1,\ldots,N$ be onedimensional projections defined with $P(X):=(0,\ldots,0,x_{k},0\ldots,0)$. Then there exists $r>0$ such that for $\left\|X_1\right\|<r$ we have $$\dim_{B}P_{k}(S(X_1))=0,\,\,\forall k\in\left\{1,\ldots,N\right\}.$$
\end{cor}
It follows from Theorem 2 andd Lemma 3.

\textbf{Remark 4.} In a case of unstable fixed point, the similar result is true. In the case of saddle, we can look the projections of orbits on stable and unstable manifolds.
\medskip

\section{Box dimension of nonhyperbolic fixed point in $\mathbb{R}^{n}$}

Now we consider the nonhyperbolic fixed points of $n$ dimensional discrete dynamical system
\begin{equation} \label{nehipsistem}
X_{n+1}=F(X_{n}).
\end{equation}

We recall Center manifold theorem in $\mathbb{R}^{n}$ .

\begin{theorem} \textbf{Center manifold},\cite{kuz}\\
Let $\mathbf{F}:\mathbb{R}^{N}\rightarrow \mathbb{R}^{N}$ be a $C^{r}$ diffeomorphism with a nonhyperbolic fixed point $\mathbf{x_0}$, with $n_0=m$. Then there is a locally defined smooth $m$-dimensional invariant manifold $W^{c}_{loc}(0)$ of the system $\mathbf{x}_{n+1}=\mathbf{G}(\mathbf{x}_{n})$ that is tangent to $T^{c}$ at $\mathbf{x}=0$. Moreover, there is a neighbourhood $U$ of $\mathbf{x}=0$ such that if $\varphi^{t}(\mathbf{x})\in U$ for all $t\geq 0$ ($t\leq 0$), then $\varphi^{t} \mathbf{x} \rightarrow W^{c}_{loc}(0)$ for $t\rightarrow +\infty$ ($t \rightarrow -\infty$).
\end{theorem}

In its eigenbasis, system ($\ref{nehipsistem}$) can be written as
\begin{eqnarray} \label{sustav11}
u&\mapsto Bu+g(u,v) \nonumber\\
v&\mapsto Cv+h(u,v)
\end{eqnarray} 
where $u\in\mathbb{R}^{m}$, $v\in\mathbb{R}^{N-m}$, matrix $B$ has order $m$ and $m$ eigenvalues on the unit circle, while the eigenvalues of $C$ are inside and/or outside it. Functions $g$ and $h$ have Taylor expansions starting with at least quadratic terms. The center manifold of the system ($\ref{sustav11}$) can be locally represented as a graph of a smooth function
$W^{c}=\{(u,v):v=V(u)\}$ where $V:\mathbb{R}^{m}\rightarrow \mathbb{R}^{N-m}$, and due to the tangent property of center manifold, $V(u)=\mathcal{O}(\left\|u\right\|^2)$.

\begin{theorem} (\textbf{Reduction principle}, \cite{kuz})\\
The system ($\ref{sustav11}$) is locally topologically equivalent near the origin to the system 
\begin{eqnarray}\label{eigens}
u&\mapsto Bu&+g(u,V(u)),\nonumber\\
v&\mapsto Cv&.
\end{eqnarray} 
\end{theorem}
The first equation in ($\ref{eigens}$) is the restriction of the system to its center manifold. The dynamics are determined by this restriction since the second equation is linear.
 
\subsection{Onedimensional center manifold in $\mathbb{R}^{n}$}

In the case of only one multiplier on the unit circle, that is $\lambda_1=\pm 1$ , the system ($\ref{sustav11}$) can be written as
\begin{eqnarray} \label{sustav12}
u&\mapsto \lambda_1 u+g(u,v) \nonumber\\
v&\mapsto Cv+h(u,v)
\end{eqnarray} 
where $u\in\mathbb{R}$, $v\in\mathbb{R}^{N-1}$, the eigenvalues of $C$ are inside and/or outside the unit circle.

\begin{theorem}
Let we have a system $\mathbf{x}\mapsto \mathbf{F}(\mathbf{x})$, $\mathbf{F}:\mathbb{R}^{N}\rightarrow \mathbb{R}^{N}$. Let $\mathbf{x_0}=0$ be a hyperbolic fixed point with $\lambda_1=\pm 1$ and $N-1$ multipliers inside or outside the unit circle. Let the center manifold be given by $u=V(u)$, where  $V:\mathbb{R}\rightarrow \mathbb{R}^{N-1}$ is of class $C^{r}$.
Let $S(u_1,v_1)=(u_{n},v_{n})$ be an orbit of the systems on the center manifold in the form 
\begin{eqnarray} \label{cenmno2}
u_{n+1}&=G(u_{n})=&\lambda_1 u_{n} + g(u_{n},V(u_{n}))\\
v_{n+1}&=V(u_{n+1})& \nonumber
\end{eqnarray}
with initial point $(u_1,v_1)$ near the origin. If the map $G$ is a $k$-nondegenerate in $x=0$, 
then there exists $r>0$ such that for $|x_1|<r$ we have $$\dim_{B}S(u_1,v_1)=1-\frac{1}{k}.$$
\end{theorem}
 
\textbf{Proof.}\\
For $u\in\mathbb{R}$ and $v\in\mathbb{R}^{N-1}$ we define a set $A=\{(u,v)\in\mathbb{R}\times\mathbb{R}^{N-1}\,:\,u\in A(u_1),\,v=0\}$, where $A(u_1)=(u_{n})_{n\in\mathbb{N}}$ is one-dimensional discrete dynamical system generated by $u_{n+1}=G(u_{n})=\Lambda_1 u_{n} + g(u_{n},V(u_{n}))$ and $\left\|(u_1,0)\right\|<r$. Then we act with the map $Q(u,z)=(u,z+G(u))$ on the set $A$ and get $$Q(A)=Q(u,0)=(u,G(u)).$$ 
So the image of the set $A$ under map $Q$ is associated to the system ($\ref{cenmno2}$) on the center manifold $Y=G(X)$. In other words, the map $Q$ 
associate the projection of an orbit on one-dimensional coordinate plane with appropriate orbit on the center manifold. 
Since the map $G$ is of class $C^{r}$, $r\geq 1$ on some neighbourhood small enough $\left\|u\right\|<\delta$, then $\left\|D G(0)\right\|\leq M$, $M>0$ ($DG$ is limited operator on that neighbourhood), and we have
$$\left\|G(u_2)-G(u_1)\right\|=\left\|D G(u^{*})\right\|\left\|u_2-u_1\right\|\leq M\left\|u_2-u_1\right\|$$
for some $u^{*}\in \left( u_1,u_2 \right)$. Therefore $G$ is a Lipschitz map for $\left\|u\right\|<\delta$.
Now it follows from Lemma 2 that $Q$ is a bilipschitz map, and 
$$\dim_{B} Q(A)= \dim_{B}A.$$

Since the map $G$ is a $k$-nondegenerate in $x=0$, from Theorem 2.2, $\cite{laho}$ it follows there
 exists $r>0$ such that for the sequence $A(u_1)=\{u_{n}\}$ defined by $u_{n+1}=G(u_{n})$ with $|x_1|<r$ we have
 $$\dim_{B}A(u_1)=1-\frac{1}{k}.$$
 We see now that $\dim_{B}A=\dim_{B}A(u_1)=1-\frac{1}{k}$.
Notice that $S(X_1,Y_1)=Q(A)$, so it follows that for $\left\|(X_1,Y_1)\right\|<r$ is
$$\dim_{B}S(X_1,Y_1)=\dim_{B} Q(A) = \dim_{B}A=1-\frac{1}{k}.\,\,\,\,\blacksquare$$
 

\subsection{Center manifold in $\mathbb{R}^{n}$}

If $\mathbf{f}:\mathbb{R}^{n}\rightarrow \mathbb{R}^{n}$ is of class $\mathcal{C}^{r}$ and $\mathbf{f}(0)=0$ and $A=D\mathbf{f}(0)=\rm{diag}[C,P,Q]$ where quadratic matrix $C$ has $c$ multipliers on the unit circle, quadratic matrix $P$ has $s$ multipliers inside the unit circle and $Q$ has $u$ multiplier outside the unit circle. Then the system can be written in a form

\begin{eqnarray} \label{cenman3}
x&\mapsto Cx+F(x,y,z) \\
y&\mapsto Pv+G(x,y,z)\nonumber\\
z&\mapsto Qz+H(x,y,z)\nonumber
\end{eqnarray}

By the center manifold theory, there exist functions $h_1$ i $h_2$ such that the system ($\ref{cenman3}$) is topologically equivalent to the system

\begin{eqnarray} \label{cenman4}
x&\mapsto Cx &+F(x,h_1(x),h_2(x)) \\
y&\mapsto Py\nonumber\\
z&\mapsto Qz\nonumber
\end{eqnarray}
for $(x,y,z)\in\mathbb{R}^{c}\times\mathbb{R}^{s}\times\mathbb{R}^{u}$ near the origin.

\begin{theorem} \textbf{Box dimension on the center manifold}\\
Let ($\ref{sustav11}$) be a system with a nonhyperbolic fixed point $\mathbf{x_0}$, matrix $B$ has $m$ multipliers on the unit circle and $C$ has $N-m$ inside the unit circle. Let the center manifold be given by $u=V(u)$, where  $V:\mathbb{R}^{m}\rightarrow \mathbb{R}^{N-m}$ is of class $C^{r}$. If an orbit $S(u_1)$ of $m$-dimensional system $u\mapsto Bu+g(u,V(u))$ has a box dimension $\dim_{B}S(u_1)=D$, then the restriction of one orbit of system ($\ref{sustav11}$) on the center manifold has s box dimension $D$. 
\end{theorem}

The proof is analogous to proof for one-dimensional manifold.

The next theorem gives the sufficient condition for nonhyperbolicity of a fixed point using the projective box dimension. 

\begin{theorem}
Let we have a system $\mathbf{x}\mapsto \mathbf{F}(\mathbf{x})$, $\mathbf{F}:\mathbb{R}^{N}\rightarrow \mathbb{R}^{N}$. Let $\mathbf{x_0}=0$ be a fixed point, and let $S(X_1)=(X_{n})_{n\in\mathbb{N}}$ be an orbit of the system $X\mapsto F(X)$ with initial point $X_1$ and $X_{n}\rightarrow 0$. If there exists $k\in\{1,\ldots,N\}$ such that $$\dim_{P_{k}}S(X_1)>0,$$ then $\mathbf{x_0}$ is a nonhyperbolic fixed point.
\end{theorem}

\textbf{Proof.}\\
If $\dim_{P_{k}}S(X_1)>0$ for some $k$, then from Lemma 3 we get that $\dim_{B}S(X_1)>0$.
The nonhyperbolicity of the fixed point $\mathbf{x_0}$ follows from Theorem 2. Namely, if $\mathbf{x_0}$ is a hyperbolic fixed point, then the box dimension of all the orbits which tend to $\mathbf{x_0}$ will be zero. Since there exists the orbit with positive box dimension we can conclude that $\mathbf{x_0}$ is a nonhyperbolic fixed point. $\blacksquare$\\

\textbf{Remark 5. }This result could be used as a tool in numerical analysis of fixed points.

\section{Box dimension of hyperbolic and nonhyperbolic singularities in $\mathbb{R}^{n}$}

In this section we use the connection of discrete and dynamical sytems via the unit-time map.
Namely, if we consider the flow of continuous system $\varphi_{t}(x)$ and put the fixed time ($t=1$), we get the map $\varphi_{1}(x)$ which generates one discrete dynamical system with the same dimension. So, the unit-time map of continuous dynamical system in $\mathbb{R}^{n}$ is a map which generates discrete dynamical system in $\mathbb{R}^{n}$ with the orbits lying on the trajectories of continuous system. In the case of hyperbolic singularity, it is obvious that every trajectory has a box dimension 1, but we will show that the unit-time map around the hyperbolic singularity has a box dimension 0.

We study the continuous dynamical system
\begin{equation} \label{k0}
\mathbf{\dot{x}}=\mathbf{F}(\mathbf{x}), \,\,\,\mathbf{x}\in\mathbb{R}^{n}.
\end{equation}
The point $\mathbf{x_0}$ is a \textbf{singularity} of a system ($\ref{k0}$) if $\mathbf{F}(\mathbf{x_0})=\mathbf{0}$. Singularity $\mathbf{x_0}$ is \textbf{hyperbolic} if there is no eigenvalue of the matrix $D\mathbf{F}(\mathbf{x_0})$ with real part equal to zero. Singularity is \textbf{nonhyperbolic} if it isn't hyperbolic. We will use the stable, unstable and center manifolds, as we did in previous sections.\\


So, we study a continuous dynamical system
\begin{equation} \label{k1}
\dot{\mathbf{x}}=\mathbf{F}(\mathbf{x})
\end{equation}
where $\mathbf{x}\in\mathbb{R}^{n}$.
The simplest way of exctracting the discrete dynamical system from the systems ($\ref{k1}$) is by using the time map $\phi_{t}(x)$. Namely, 
we fix $t_0>0$ and observe the system on $X$ which is generated by the iteration of the $t_0$-time map $\phi_{t_0}$. We simply put $t_0=1$, and get the discrete system generated by the unit-time map
\begin{equation} \label{k2}
\mathbf{x} \mapsto \phi_{1}(\mathbf{x}).
\end{equation} 
It is easy to show that the isolated fixed point of the system ($\ref{k2}$) correspond to isolated singularities of the system ($\ref{k1}$). In order to study that correspodance between the hyperbolicity and nonhyperbolicity of appropriate points, we need to find the connection between eigenvalues of $DF(\mathbf{x_0})$ i $D\phi_{1}(\mathbf{x_0})$. It is easily seen that $D\phi_1(\mathbf{x_0})=e^{A}$ where $A=DF(\mathbf{x_0})$. It means that Jacobi matrices are connected in a way that $\mathbf{x_0}$ is hyperbolic (nonhyperbolic) singularity of the system $(\ref{k1})$ if and only if $x_0$ is a hyperbolic (nonhyperbolic) fixed point of the map ($\ref{k2}$).

\begin{theorem} \textbf{Stable hyperbolic singularity in $\mathbb{R}^{n}$}\\
Let $\mathbf{x_0}=0$ be a hyperbolic singularity of a system ($\ref{k1}$) with all the eigenvalues with negative real part, and let $\phi_1(\mathbf{x})$ be an appropriate unit-time map. Then there exists $r>0$ such that for the sequence $S(\mathbf{x}_1)=(\mathbf{x}_{k})_{k\in\mathbb{N}}$ defined with $\mathbf{x}_{k+1}=\phi_1(\mathbf{x}_{k})$, $\left\|\mathbf{x}_{0}\right\|<r$ we have $\dim_{B} S(\mathbf{x}_1)= 0$.
\end{theorem}
\textbf{Proof.} From the Taylor series for 
$$\phi_1(\mathbf{x})=e^{A}\mathbf{x}+g^{k}(\mathbf{x})+\mathcal{O}(\left\|x\right\|^{k+1})$$ where $A=DF(\mathbf{x}_0)$
we have that $D\phi_1(\mathbf{x}_0)=e^{A}$, so it follows that tha stable hyperbolic singularity of the system ($\ref{k1}$) is also a stable hyperbolic fixed point of appropriate unit-time map. From Theorem 2 we have $\dim_{B}S(x_1)=0$. $\blacksquare$\\

\textbf{Remark 6} Analogously for the unstable case $Re(\lambda)>0$. In the case of hyperbolic saddle, we apply the theorem on the stable and unstable manifolds.\\

It is clear that the form of every one-dimensional continuous system with nonhyperbolic singularity $x_0=0$ is
\begin{equation} \label{k04}
\dot{x}=f^{(m)}(x)+\mathcal{O}(\left|x\right|^{m+1}).
\end{equation}
for some $m>1$. The next lemma gives the form of the appropriate unit-time map. 

\begin{lemma}
Let continuous dynamical system $$\dot{x}=f^{(m)}(x)+f^{(m+1)}(x)+\ldots,\,\,\,x\in\mathbb{R},$$ with $m>1$, has a nonhyperbolic singularity in $x_0=0$.
The the unit-time map has a form $$\phi_1(x)=x+f^{(m)}(x)+\mathcal{O}(\left|x\right|^{m+1}).$$
\end{lemma}


At this point, we study the orbits of the discrete dynamical system generated by the unit-time map, that is, the sequences
 $S(x_1)=(x_{n})_{n\geq 1}$ defined with
\begin{equation} \label{stabtok}
x_{n+1}=\left\{\begin{array}{ll}
\phi_1(x_{n}), & \textrm{if } x_0 \textrm{ is stable};\\
\phi^{-1}(x_{n}), & \textrm{if } x_0 \textrm{ is unstable}.
\end{array} \right.
\end{equation}

\begin{theorem} \textbf{Nonhyperbolic singularity in $\mathbb{R}$}\\
Let $x_0=0$ be a nonhyperbolic singularity of the system
($\ref{k04}$), and let $\phi_1(x)$ be associated unit-time map. 
Then there exists $r_1>0$ such that for the sequence $S(x_1)=(x_{n})_{n\geq 1}$ defined by ($\ref{stabtok}$), $x_1\in(0,r_1)$ we have
$$\dim_{B}S(x_1)=1-\frac{1}{m}.$$
\end{theorem}
\textbf{Proof.}\\
We apply Lemma 4, Theorem 1 and Colorally 2. $\blacksquare$\\

We observe the continuous planar system with the eigenvalues: $\lambda_1=0$ and $\lambda_2$ is a negative real number. 
Then the form of the system is
\begin{eqnarray} \label{ks2}
\dot{x}&=&f(x,y)\nonumber\\
\dot{y}&=&\lambda_2{y}+g(x,y),
\end{eqnarray}
By Center manifold theorem we get that the restriction of the system ($\ref{ks2}$) on the one-dimensional local central $C^{r}$ manifold $y=h(x)$ is of a form
\begin{eqnarray}\label{cm1}
\dot{x}&=&f(x,h(x))\nonumber\\
\dot{y}&=&\lambda_2{y}
\end{eqnarray}
for all $x\in \mathbb{R}$ for which $|x|<\delta$.

In the new theorem we consider the two-dimensional sequences $S(x_1,y_1)=\{(x_{n},y_{n})\}_{n\geq 1}$ defined by 
\begin{eqnarray} \label{stabtok2}
x_{n+1}&=&\left\{\begin{array}{ll}
\phi_1(x_{n}), & \textrm{if } x_0 \textrm{ is stable}\\
\phi^{-1}(x_{n}), & \textrm{if } x_0 \textrm{ is unstable}
\end{array} \right. \nonumber\\
y_{n+1}&=&h(x_{n+1}).
\end{eqnarray}
This is the two-dimensional discrete dynamical system generated with the unit-time map of the trajectory on the center manifold.

\begin{theorem} \textbf{Nonhyperbolic singularity in $\mathbb{R}^2$}\\
Let a system ($\ref{ks2}$) has a restriction on the local center manifold given by ($\ref{cm1}$).
Let $\phi_1(x)$ be an unit-time map of the system ($\ref{cm1}$), and assume that the two-dimensional sequence $S(x_1,y_1)=(x_{n},y_{n})_{n\geq 1}$ is defined by 
($\ref{stabtok2}$) with the initial point $(x_1,y_1)$. If the map $\phi_1$ is a $k$-nondegenerate in a point $x_0=0$, then there exists
$r>0$ such that for $\left\|(x_1,y_1)\right\|<r$ we have $$\dim_{B}S(x_1,y_1)=1-\frac{1}{k}.$$
\end{theorem}
\textbf{Proof.}\\
As in the proof of Theorem 12, we use Lemma 4 and the fact that the map $y=h(x)$ is of class $C^{r}$ and Lipschitz on the neighbourhood of $x_0$. Then the claim easily foolows form Theorem 1 and Colorally 2. $\blacksquare$\\


\begin{cor}
Let $\mathbf{x_0}=0$ be a singularity of a system ($\ref{k1}$), and let $\phi_1(\mathbf{x})$ be an appropriate unit-time map. 
Let $S(X_1)=(X_{n})_{n\in \mathbb{N}}$ be an orbit of the system $X\mapsto \phi_1(X)$ with initial point $X_1$ and $X_{n}\rightarrow 0$. If there exists $k\in\{1,\ldots,N\}$ such that $$\dim_{P_{k}}S(X_1)>0,$$ then $\mathbf{x_0}$ is a nonhyperbolic singularity.
\end{cor}

This claim easily follows from Theorem 10.










\section*{Acknowledgment}
This article was supported by the Croatian Science Foundation Project IP-2014-09-2285.

\end{document}